\documentclass[a4paper,12pt]{article}
\usepackage{amssymb,amsmath}

\newcommand{\R}{\mathbb{R}}

\newcommand{\N}{\mathbb{N}}
\newcommand{\ep}{\varepsilon}
\newcommand{\x}{\mathbf{x}}
\newcommand{\Z}{\mathbb{Z}}
\renewcommand{\mod}[1]{\bmod #1}
\renewcommand{\b}[1]{\mathbf{#1}}

\newtheorem{lemma}{Lemma}
\newtheorem{theorem}{Theorem}

\begin{document}

\title{Powerfree Values of Polynomials}
\author{D.R. Heath-Brown\\Mathematical Institute, Oxford}
\date{}
\maketitle

\section{Introduction}\label{intro}

Let $f(X)\in\Z[X]$ be an irreducible polynomial of degree $d$.  It is
conjectured that, for any integer $k\ge 2$, the polynomial $f(n)$
takes infinitely many $k$-th power free values, providing that $f$
satisfies the obviously necessary congruence conditions. Thus for
every prime $p$ we need to assume that there is at least one integer
$n_p$ for which $p^k\nmid f(n_p)$. This problem appears to become
harder as the degree $d$ increases, but easier as $k$ increases.  Thus
in 1933 Ricci \cite{ricci} handled the case $k\ge d$, and even proved an
asymptotic formula
\begin{equation}\label{as}
N_{f,k}(x)\sim A(f,k)x\;\;\;(x\rightarrow\infty)
\end{equation}
where
\[N_{f,k}(x):=\#\{n\in\N:\, n\le x,\, f(n) \mbox{ $k$-free}\}.\]
Here the constant $C(f,k)$ is given as
\[C(f,k):=\prod_p(1-\rho_f(p^k)p^{-k})\]
where
\[\rho_f(d):=\#\{n\mod{d}:\, d\mid f(n)\}.\]
Further progress was made twenty years later by Erd\H{o}s 
\cite{erd} who showed that one
could obtain $k$-free values for $k=d-1$, as soon as $d\ge 3$.  For
such $k$ the asymptotic formula (\ref{as}) was later obtained by
Hooley \cite{hoo}.  The next development was due to Nair \cite{nair}
who established (\ref{as}) for $k\ge (\sqrt{2}-\frac{1}{2})d$.  In
particular Nair's result shows that $k=d-2$ is admissible for $d\ge
24$.  The author \cite[Theorem 16]{HB1}
then showed how the ``determinant method'' could be
applied to the problem, and demonstrated that the asymptotic formula
remained valid for $k\ge (3d+2)/4$, so that one may take $k=d-2$
providing only that $d\ge 10$.  Indeed using methods of Salberger 
(to appear) one can replace these inequalities by $k\ge (3d+1)/4$ and
$d\ge 9$ respectively.

In this paper we show that further progress is possible for irreducible
polynomials of the form $f(X)=X^d+c$.  For these we establish the following
result.
\begin{theorem}\label{mainthm}
Let $f(X)=X^d+c\in\Z[X]$ be an irreducible polynomial, and suppose
that $k\ge (5d+3)/9$.  Then there is a constant $\delta(d)$ such that
\[N_{f,k}(x)= C(f,k)x+O(x^{1-\delta(d)}).\]
The implied constant may depend on $f$ and $k$.
\end{theorem}
For comparison with the earlier results we point out that this will
allow $k=d-2$ as soon as $d\ge 6$. The result of Erd\H{o}s handles the
case of cubic polynomials taking square-free values, and the most
interesting open question then concerns quartic polynomials taking
square-free values.  We would therefore like to handle $k=d-2$ for
$d=4$, and one can track our progress towards this goal through the
historical discussion above.

There is a related question concerning powerfree values of $f$ at
prime arguments.  Here there is a natural condition that for every
prime $p$ there should be an integer $n_p$, coprime to $p$, and such
that $p^k\nmid f(n_p)$. With this in mind one defines
\[N'_{f,k}(x):=\#\{p \mbox{ prime}:\, p\le x,\, f(p) \mbox{ $k$-free}\}\]
and
\[C'(f,k):=\prod_p(1-\rho_f'(p^k)\phi(p^k)^{-1})\]
where
\[\rho_f'(d):=\#\{n\mod{d}:\, {\rm g.c.d.}(n,d)=1,\, d\mid f(n)\}.\]
The corresponding asymptotic formula
\[N'_{f,k}(x)\sim C'(f,k)\pi(x)\;\;\;(x\rightarrow\infty)\]
has been proved for $k=d$ by Uchiyama \cite{Uch}, by a method that
also handles the case $k>d$. However it remains an open problem to
establish this in the case $k=d-1$ considered for the previous problem
by Erd\H{o}s and Hooley.  None the less, important progress has been
made by Helfgott \cite{helf1} and \cite{helf2}, showing in particular
that the asymptotic formula holds in the case $k=2$ and $d=3$.

Our methods are sufficiently robust that they apply immediately to
powerfree values of $f(p)$.  We have the following result.
\begin{theorem}\label{primethm}
Let $f(X)=X^d+c\in\Z[X]$ be an irreducible polynomial, and suppose
that $k\ge (5d+3)/9$.  Suppose that for every
prime $p$ there is an integer $n_p$, coprime to $p$, and such
that $p^k\nmid f(n_p)$. Then for any fixed $A>0$ we have
\[N'_{f,k}(x)= C'(f,k)\pi(x)+O_A(x(\log x)^{-A}).\]
In particular this holds for $k=d-1$ and every $d\ge 3$.
\end{theorem}

The preliminary manoeuvres for these problems are straightforward. We
shall fix the polynomial $f$ (and hence also $d$) throughout, so that
all order constants may depend tacitly on $f$ and $d$.  The key fact
we shall use is that
\[\sum_{b^k|f(n)}\mu(b)=\left\{\begin{array}{cc}1, & f(n)\;\mbox{is
$k$-free},\\ 0, & \mbox{otherwise}.\end{array}\right.\]
It follows that
\[N_{f,k}(x)=\sum_b \mu(b)N(b,x)\]
with
\[N(b,x)=\#\{n\le x: b^k|f(n)\},\]
and similarly that
\[N'_{f,k}(x)=\sum_b \mu(b)N'(b,x)\]
with
\[N'(b,x)=\#\{p\le x: b^k|f(p)\}.\]
Clearly $N'(b,x)=N(b,x)=0$ for $b\gg x^{d/k}$.  If we denote the
solutions to $f(n)\equiv 0\mod{b^k}$ by $n_1,\ldots,n_r$, where
$r=\rho_f(b^k)$, then
\begin{eqnarray*}
N(b,x)&=&\sum_{i\le r}\#\{n\le x: n\equiv n_i\mod{b^k}\}\\
&=&\sum_{i\le r}(xb^{-k}+O(1))\\
&=&xb^{-k}\rho_f(b^k)+O(\rho_f(b^k)),
\end{eqnarray*}
and similarly, providing that $b\le (\log x)^{2A}$ we have
\begin{eqnarray*}
N'(b,x)&=&\sum_{i\le r}\#\{p\le x: p\equiv n_i\mod{b^k}\}\\
&=&\sum_{i\le r,\, (n_i,b)=1}\pi(x;b^k,n_i)\\
&=&\frac{\pi(x)}{\phi(b^k)}\rho_f'(b^k)+O_A(\rho_f(b^k)x(\log x)^{-4A}),
\end{eqnarray*}
by the Siegel-Walfisz Theorem.  Now, for any $\xi>0$  we find that
\[\sum_{b\le \xi}\mu(b)N(b,x)=
x\sum_{b\le\xi}\frac{\mu(b)\rho_f(b^k)}{b^k}+O(\sum_{b\le\xi}\rho_f(b^k)).\]
The function $\rho_f$ is multiplicative, with $\rho(p^k)\ll 1$,
whence \begin{equation}\label{rhoest}
\rho(b^k)\ll_{\ep} b^{\ep}
\end{equation}
for any $\ep>0$ and any square-free $b$.  If $k\ge 2$ 
it follows on taking $\ep=1/2$ that
\[\sum_{b\le\xi}\frac{\mu(b)\rho_f(b^k)}{b^k}=
\sum_{b=1}^{\infty}\frac{\mu(b)\rho_f(b^k)}{b^k}+O(\sum_{b>\xi}b^{1/2-k})
=C(f,k)+O(\xi^{-1/2})\]
and
\[\sum_{b\le\xi}\rho_f(b^k)\ll \xi^{3/2}.\]
In particular if we set $\xi=x^{1/2}$ we see that
\[\sum_{b\le \xi}\mu(b)N(b,x)=C(f,k)x+O(x^{3/4}).\]
In precisely the same way, if we take $\xi=(\log x)^{2A}$, then
\[\sum_{b\le \xi}\mu(b)N'(b,x)=C'(f,k)\pi(x)+O_A(x(\log x)^{-A}).\]

We now consider the range $\xi<b\le x^{1-\eta}$, where $\eta$ is a
small positive constant.  Here we have
\[\sum_{\xi<b\le x^{1-\eta}}\mu(b)N(b,x)\ll
\sum_{\xi<b\le x^{1-\eta}}N(b,x)\ll
\sum_{\xi<b\le
  x^{1-\eta}}\left(\frac{x}{b^k}+O(1)\right)\rho_f(b^k).\]
If we use the bound (\ref{rhoest}) with $\ep=\tfrac12 \eta\le 
\tfrac12$ this yields
\[\sum_{\xi<b\le x^{1-\eta}}\mu(b)N(b,x)\ll x\xi^{-1/2}+x^{1-\eta/2}
\ll x^{1-\eta/2}.\]
This bound is satisfactory for Theorem \ref{mainthm}.
Since $N'(b,x)\le N(b,x)$ we will get exactly the same bound in the estimation
of $N'_{f,k}(x)$, and again this is satisfactory for Theorem \ref{primethm}.

To complete the proof of the two theorems it will now be enough to
show that
\[\sum_{x^{1-\eta}<b\ll x^{d/k}}N(b,x)\ll x^{1-\delta}\]
for some $\delta>0$, providing that $\eta$ is small enough.
By a suitable dyadic subdivision we then see that it will suffice to
establish the following estimate.
\begin{lemma}\label{mainlemma}
Let $f(X)=X^d+c\in\Z[X]$ be an irreducible polynomial. For any
$N,A,B\in\N$ define
\begin{eqnarray*}
\lefteqn{F(N;A,B):=}\\
&&\#\{(n,a,b)\in\N^3: f(n)=ab^k , N<n\le 2N,\,A<a\le 2A,\,
B<b\le 2B\}.
\end{eqnarray*}
Then if $(5d+3)/9\le k \le d-1$ there is a constant $\delta$ depending 
on $d$ such that
\[F(N;A,B)\ll_f N^{1-\delta}\]
for $B\ge N^{1-\delta}$.
\end{lemma}

We have now reduced our problem to one of counting solutions to a
Diophantine equation $f(n)=ab^k$, inside a suitable box.  A general
procedure for such questions is provided by the ``determinant method''
developed in the author's paper \cite{ann}.  The efficiency of the
method depends on the dimension of the associated algebraic variety.
For $f(n)=ab^k$ we are counting integer points on an affine surface.
Thus far we have made no use of the special shape of the polynomial
$f$, but if we observe that $f(n)=n^d+O(1)$ we see that $(n,a,b)$ lies
close to the weighted projective curve $X_0^d=X_1X_2^k$, where
$X_0$ and $X_2$ are given weight 1, and $X_1$ has weight $d-k$.  Thus
the particular form of the polynomial $f$ allows us to consider points
close to a curve, rather that points on a surface. Reducing the
dimension in this way is the key to our saving.  The procedure is
discussed in more detail in the author's work \cite{3powers}, to which
the interested reader should be directed.

\section{The Determinant Method}

Since $f(n)=n^d+O(1)$ we will have 
\begin{equation}\label{Ab}
N^dB^{-k}\ll A\ll N^dB^{-k}
\end{equation}
for large $N$. Moreover, since $a\ge 1$ we may assume that $B^k\ll
N^d$, and indeed we shall assume that
\begin{equation}\label{Br}
N^{1-\eta}\ll B\ll N^{d/k}
\end{equation}
for some positive constant $\eta$.
We will choose a parameter $K\ge 1$ having
\begin{equation}\label{kr}
1\ll\frac{\log K}{\log N}\ll 1,
\end{equation}
and divide the 
available range for $n/b$ into $O(K)$ subintervals 
\[I=(m_0N/BK,(m_0+1)N/BK]\]
with endpoints defined by integers $m_0$ in the range
\begin{equation}\label{m0}
K\ll m_0\ll K.  
\end{equation}
We use $F_I(N;A,B)$
to denote the corresponding contribution to $F(N;A,B)$.  Since
$f(n)=n^d+O(1)$ we have $n^d=ab^k+O(1)$ and
\[(n/b)^d=a/b^{d-k}+O(B^{-d}).\]
It will be convenient to put $k=d-j$ so that 
\begin{equation}\label{M}
(n/b)^d=a/b^{j}+O(B^{-d}).
\end{equation}

We now begin the determinant method by listing the points
$(n_r,a_r,b_r)$ contributing to $F_I(N;A,B)$.  Thus the index $r$ runs
from 1 to 
\[R:=F_I(N;A,B).\]
We choose an integer parameter $D\ge 1$ and 
consider the monomials
\[m(n,a,b)=n^ua^vb^w\]
for which $u+jv+w=D$.  Thus we may consider $D$
as the weighted degree of the monomial, where the variables $(n,a,b)$ are
given weights $(1,j,1)$.  The number of such monomials will be
\begin{equation}\label{H}
H:=\sum_{v\le D/j}(D-jv+1)=\frac{D^2}{2j}+O(D)
\end{equation}
and we label them as $m_1(n,a,b),\ldots,m_H(n,a,b)$.  We now proceed
to consider the $R\times H$ matrix $M$ say, whose $(r,h)$ entry is
$m_h(n_r,a_r,b_r)$.  The strategy of the determinant method is to show
that $M$ has rank strictly less than $H$, if the parameters $K$ and
$D$ are suitably chosen.  If this can be achieved, there will be a
non-zero integer vector $\b{c}$ such that $M\b{c}=\b{0}$. This vector
will depend on the interval $I$, that is to say it will depend on
$m_0$. It
provides the coefficients of a weighted homogeneous polynomial 
\[C_I(n,a,b)=\sum_{h}c_hm_h(n,a,b)\]
such that 
\begin{equation}\label{aux}
C_I(n_r,a_r,b_r)=0,\;\;\;(r\le R).
\end{equation}

If $R<H$ the matrix $M$ automatically has rank less than $H$.
Otherwise it suffices to show that any $H\times H$ sub-determinant
vanishes, and it will be enough to consider the determinant formed from
the first $H$ rows of $M$, which we shall denote by $\Delta$.  Clearly
$\Delta$ is an integer, and our strategy is to show that $|\Delta|<1$
so that $\Delta$ must vanish.

We proceed to divide the $r$-th row of $\Delta$ by $b_r^DB^{-D}$ for each
$r\le D$, and similarly to divide
the column corresponding to the monomial $n^ua^vb^w$ by
$N^uA^vB^{w}$.  Since 
\[n^ua^vb^w
=\left(\frac{b}{B}\right)^D\left(\frac{nB}{bN}\right)^u
\left(\frac{aB^j}{b^jA}\right)^vN^uA^vB^w\]
for $u+jv+w=D$, this produces a new determinant $\Delta_1$ whose
entries are of the form $m_h(nB/bN,aB^j/b^jA,1)$.  Moreover we have
\begin{equation}\label{D1}
|\Delta|=|\Delta_1|\prod_{r\le H}(b_r/B)^D\prod_{u,v,w}N^uA^vB^w\le 
2^{HD}P|\Delta_1|,
\end{equation}
where
\[P=\prod_{u+jv+w=D}N^uA^vB^w.\]
If we write $B=N^{\beta}$ then we have $\log A=(d-k\beta)\log N+O(1)$,
by (\ref{Ab}).  It follows that
\begin{eqnarray}\label{P}
\log P
&=&(\log N)\sum_{u+jv+w=D}(u+v(d-k\beta)+w\beta)+O_D(1)\nonumber\\
&=&(\log N)\left\{\frac{D^3}{6j}(1+(d-k\beta)j^{-1}+\beta)+O(D^2)\right\}
+O_D(1).
\end{eqnarray}

We now write 
\[\frac{n_rB}{b_rN}=\frac{m_0}{K}+s_r,\;\;\;\mbox{and}\;\;\;
\frac{a_rB^j}{b_r^jA}=\frac{N^d}{AB^k}\left(\frac{m_0}{K}+s_r\right)^d+t_r.\]
Since $n_r/b_r\in(m_0N/BK\,,\,(m_0+1)N/BK]$ it follows that
\[s_r\ll K^{-1}.\]
Moreover (\ref{Ab}) and (\ref{M}) yield
\[\frac{a_rB^j}{b_r^jA}=
\frac{N^d}{AB^k}\left(\frac{n_rB}{b_rN}\right)^d+O(N^{-d}),\]
and hence
\[t_r\ll N^{-d}.\]
Thus the $(r,h)$ entry of $\Delta_1$ will be a polynomial 
\[f_h(s_r,t_r)=
(m_0K^{-1}+s_r)^u\left(N^dA^{-1}B^{-k}(m_0K^{-1}+s_r)^d+t_r\right)^v.\]
Clearly $f_h$ may depend on 
$h,m_0,K,D$ and $d$, but it is independent of $r$.
Moreover the degree of $f_h$ will be at most $dD$.  It follows from
(\ref{Ab}) and (\ref{m0}) that $N^dA^{-1}B^{-k}\ll 1$ and
$m_0K^{-1}\ll 1$, whence we have the bound
$||f_h||\ll_{D}1$ for the height of $f_h$.

In order to estimate the size of $\Delta_1$ we will use Lemma 3 of the
author's work \cite{3powers}.  For each of the monomials $s^ut^v$ we write
\[||s^ut^v||=K^{-u}N^{-dv},\]
and we list them in order as
$m_1,\ldots,m_H$ with $||m_1||\ge ||m_2||\ge\ldots$.  Then according
to \cite[Lemma 3]{3powers} we have
\begin{equation}\label{D2}
\Delta_1\ll_{D}(\max ||f_h||)^H\prod_{h=1}^H||m_h||\ll_{D}\prod_{h=1}^H||m_h||.
\end{equation}
To proceed further we shall write $K=N^{\kappa}$, and note that
$1\ll\kappa\ll 1$, by (\ref{kr}).
If we now write $m(\lambda)$, say, for the number of
monomials $m_r=s^ut^v$ with $||m_r||\ge N^{-\lambda}$, then
\[m(\lambda)=\#\{(u,v)\in\Z^2:\, u,v\ge 0,\, \kappa u+dv\le\lambda\}
=\frac{\lambda^2}{2\kappa d}+O(\lambda)+O(1).\]
If $||m_H||=N^{-\lambda_0}$ then $m(\lambda_0)\ge H$, while for any
$\ep>0$ we will have
\[m(\lambda_0-\ep)\le H-1.\]
We may therefore deduce that
\[\lambda_0=\sqrt{2\kappa dH}+O(1).\]
We then find that
\[\prod_{h=1}^H||m_h||=N^{-\mu}\]
with
\begin{eqnarray*}
\mu&=&\sum_{\kappa u+dv\le\lambda_0}(\kappa u+dv)+O(\lambda_0^2)+O(1)\nonumber\\
&=&\frac{\lambda_0^3}{3\kappa d}+O(\lambda_0^2)+O(1)\nonumber\\
&=&\frac{2^{3/2}}{3}(\kappa d)^{1/2}H^{3/2}+O(H).
\end{eqnarray*}
In view of (\ref{H}), 
(\ref{D1}), (\ref{P}) and (\ref{D2}) we may now conclude that
\[\frac{\log|\Delta|}{\log N}\le
\frac{D^3}{6j}(1+(d-k\beta)j^{-1}+\beta)
-\frac{2^{3/2}}{3}(\kappa d)^{1/2}H^{3/2}+O_D((\log N)^{-1})+O(D^2).\]
Thus (\ref{H}) yields
\begin{eqnarray*}
\frac{\log|\Delta|}{D^3\log N}&\le&
\frac{1}{6j}(1+(d-k\beta)j^{-1}+\beta)
-\frac{2^{3/2}}{3}(\kappa d)^{1/2}(2j)^{-3/2}\\
&&\hspace{1cm}\mbox{}+O_D((\log N)^{-1})+O(D^{-1}).
\end{eqnarray*}
We therefore choose
\begin{equation}\label{kappad}
\kappa=\frac{j}{4d}\left(1+\frac{d-k\beta}{j}+\beta\right)^2+\eta,
\end{equation}
with the same small constant $\eta$ as in (\ref{Br}).  Then
(\ref{kr}) will be satisfied, and we will have
\[\frac{\log|\Delta|}{D^3\log N}<0\]
providing that we first choose $D=D(f,d,\eta)$ sufficiently large, 
and then ensure that $N$
is sufficiently large in terms of $f,d$ and $\eta$.  

We therefore deduce that $\Delta=0$ when $K=N^{\kappa}$.
With this choice the matrix $M$ introduced at the beginning of the
section will have rank strictly less than $H$, so that all solutions
$(n_r,a_r,b_r)$ counted by $F_I(N;A,B)$ satisfy the auxiliary
equation (\ref{aux}).

\section{Completion of the Proof}

We now complete our estimation of $F_I(N;A,B)$ by considering how many
triples $(n,a,b)$ can satisfy both the original equation $f(n)=ab^k$ and
the additional equation (\ref{aux}).  The procedure here will follow
precisely that used in the author's paper \cite[\S 5.3]{HB1}.  Since $C_I$ is
homogeneous with exponent weights $(1,j,1)$ any factor would have to
be similarly weighted-homogeneous.  It follows in particular that
$C_I(x,y,z)$ cannot have a factor in common with $f(x)-yz^k$.  As in
\cite[pages 84 and 85]{HB1} we find that either
\begin{equation}\label{alt}
F_I(N,A,B)\ll_{\ep}(1+N/B)N^{\ep}
\end{equation}
or that there is an irreducible polynomial $G_I(X,Y)\in\Z[X,Y]$, with
degree bounded in terms of $d$ and $\ep$, but at least $d$, such that 
\begin{equation}\label{Geq}
G_I(n,b)=0
\end{equation}
for every triple $(n,a,b)$
counted by $F_I(N,A,B)$.

For a given interval $I$ we will have 
\[n/b\in I=(m_0N/BK,(m_0+1)N/BK]\]
It therefore follows that 
\begin{equation}\label{reg}
\left|n-\frac{m_0N}{BK}b\right|\le \frac{2N}{K},\;\;\; B<b\le 2B.
\end{equation}
It will be convenient to define a linear mapping
$T:\R^2\rightarrow\R^2$ by
\[T\x:=\left(\begin{array}{c} K(2N)^{-1}x_1-(2B)^{-1}m_0x_2 \\
(2B)^{-1}x_2\end{array}\right)\]
and to consider the lattice
\[\Lambda=\{T\x:\, \x\in\Z^2\}\] 
of determinant $K(4NB)^{-1}$. Then if $\x=(n,b)$
satisfies (\ref{reg}) we produce a point
$T\x=(\alpha_1,\alpha_2)\in\Lambda$ falling in the square 
\[S=\{(\alpha_1,\alpha_2): \max(|\alpha_1|,|\alpha_2|)\le 1\}. \]
Let $\b{g}^{(1)}$ be the shortest non-zero vector in the lattice and
$\b{g}^{(2)}$ the shortest vector not parallel to $\b{g}_1$.  These vectors
will form a basis for $\Lambda$.  Moreover we have 
$\lambda_1\b{g}^{(1)}+\lambda_2\b{g}^{(2)}\in S$ only when $|\lambda_1|\ll
|\b{g}^{(1)}|^{-1}$ and $|\lambda_2|\ll|\b{g}^{(2)}|^{-1}$. These  
constraints may be written in the form
$|\lambda_i|\le L_i$, for appropriate bounds $L_1,L_2$.
Since $|\b{g}^{(2)}|\ge|\b{g}^{(1)}|$ and 
$|\b{g}^{(1)}|.|\b{g}^{(2)}|\ll {\rm det}(\Lambda)\ll K(NB)^{-1}$ we will have
$L_1\gg L_2$ and $L_1L_2\gg NBK^{-1}$.  We now write
$\b{h}^{(i)}=T^{-1}\b{g}^{(i)}$ 
for  $i=1,2$.  These vectors will then be a basis for $\Z^2$, and if
$\b{x}=\lambda_1\b{h}^{(1)}+\lambda_2\b{h}^{(2)}$ is in the region
(\ref{reg}) then
we will have $|\lambda_i|\le L_i$ for $i=1,2$.  This allows us to make
a change of basis, replacing $(x_1,x_2)$ by $(\lambda_1,\lambda_2)$
so that our constraints on $n,b$ are replaced by the conditions
$|\lambda_i|\le L_i$.

We therefore proceed to substitute $\lambda_1,\lambda_2$ for $n,b$ in
(\ref{Geq}). We may then use the bound of Bombieri and Pila
\cite[Theorem 5]{BP}
to show that the number of possible pairs $\lambda_1,\lambda_2$ is
$\ll_{\ep}\max(L_1,L_2)^{1/d+\ep}\ll_{\ep}L_1^{1/d+\ep}$, since the
degree of $G_I$ is at least $d$.  Thus
\[F_I(N,A,B)\ll_{\ep}L_1^{1/d+\ep}.\]
The number $L_1$ depends on the interval $I$, which is determined by
$m_0$.  We therefore write $L_1=L_1(m_0)$ accordingly.  In view of the
alternative (\ref{alt}) we then see that
\begin{equation}\label{FEST}
F(N,A,B)\ll_{\ep} K(1+N/B)N^{\ep}+\sum_{K\ll m_0\ll K}L_1(m_0)^{1/d+\ep},
\end{equation}
the range for $m_0$ being given by (\ref{m0}).

We proceed to investigate the number of choices for $m_0$ which produce
a value $L_1(m_0)$ lying in a given dyadic interval $(L,2L]$ say.  
In the notation above,
if $(n,b)=(x_1,x_2)$ corresponds to $\b{g}^{(1)}$ then 
\[L_1\left(x_1-\frac{m_0N}{BK}x_2\right)\ll \frac{N}{K}\]
and $L_1x_2\ll B$.  Moreover we will have ${\rm g.c.d.}(x_1,x_2)=1$.
Thus the number of intervals $I$ for which $L<
L_1\le 2L$ is at most the number of triples $(x_1,x_2,m_0)\in\Z^3$
with ${\rm g.c.d.}(x_1,x_2)=1$, for which
\begin{equation}\label{cs}
L\left(x_1-\frac{m_0N}{BK}x_2\right)\ll \frac{N}{K},\;\;\;
Lx_2\ll B,\;\;\;\mbox{and}\;\;\;K\ll m_0\ll K.
\end{equation}
We proceed to consider whether the value $x_2=0$ can occur.  If
$x_2=0$ the first of the conditions above would yield
$Lx_1\ll N/K$.  However we cannot have $x_1=x_2=0$, so that we must
have $L\ll N/K$ whenever $x_2=0$.  We now recall that 
$L_1\gg L_2$ and that $L_1L_2\gg NBK^{-1}$, whence 
\begin{equation}\label{LC}
L^2\gg NBK^{-1}.
\end{equation}
It follows that if $x_2=0$ then $(N/K)^2\gg L^2\gg NBK^{-1}$ and hence
that $BK\ll N$.
However, since $K=N^{\kappa}$ with $\kappa$ given by (\ref{kappad}),
we see from (\ref{Br}) that $BK/N$ tends to infinity with $N$,
which ensures that the case $x_2=0$ cannot arise.

We now see in particular that the second condition of (\ref{cs}) 
yields $L\ll B$. If we rewrite the first of the 
conditions (\ref{cs}) to say that
\[m_0x_2=N^{-1}BKx_1+O(BL^{-1})\]
we then
see that each choice for $x_1$ restricts the product $m_0x_2$ to an
interval of length $\ll B/L$, with $B/L\gg 1$.  Moreover $m_0x_2$ is
never zero.  Thus a divisor function estimate shows that there are
$O_{\ep}(N^{\ep}BL^{-1})$ possible pairs $(x_2,m_0)$ for each value of
$x_1$.  The conditions (\ref{cs}) show that $x_1\ll N/L$, so that
$x_1$ takes $O(1+N/L)$ values.  This allows us to conclude that the
number of integers for $m_0$ which produce a value $L_1(m_0)$ in the range
$L< L_1\le 2L$ is $O_{\ep}((1+N/L)N^{\ep}BL^{-1})$.

We can now feed this information into (\ref{FEST}), using a dyadic
subdivision for the values of $L_1(m_0)$ to obtain 
\[F(N,A,B)\ll_{\ep} K(1+N/B)N^{\ep}+\sum_{L}L^{1/d+\ep}(1+N/L)N^{\ep}BL^{-1},\]
in which $L$ runs over powers of 2, subject to the condition $L\gg
(NBK^{-1})^{1/2}$ given by (\ref{LC}).  It then follows that
\[F(N,A,B)\ll_{\ep} K(1+N/B)N^{\ep}+L_0^{1/d+\ep}(1+N/L_0)N^{\ep}BL_0^{-1},\]
where $L_0:=\max\{1\,,\,(NBK^{-1})^{1/2}\}$.  
On taking $\ep=\eta$ we deduce from (\ref{Br}) that 
\[F(N,A,B)\ll_{\eta} N^{2\eta}\{K+L_0^{1/d+\eta}(1+N/L_0)BL_0^{-1}\}.\]

We proceed to analyse our estimate for $F(N,A,B)$ by defining
\[\rho(t)=\frac{j}{4d}\left(1+\frac{d-kt}{j}+t\right)^2\]
and $q(t)=\rho(t)+1-t$. Then 
\[q'(t)=
\frac{j}{d}\left(1+\frac{d-kt}{j}+t\right)\left(1-\frac{k}{j}\right)-1.\]
This is clearly negative if $k\ge j$ and $0\le t\le d/k$.  Hence if
$k\ge d/2$ we have 
\[q(t)\ge q(d/k)=\frac{j^3}{4dk^2}\ge 0\]
for $0\le t\le d/k$.  It therefore follows that $KN\le B$, and hence
that $L_0\le N$ for the relevant range of $B$.  Our estimate now
simplifies to give
\[F(N,A,B)\ll_{\eta} N^{2\eta}\{K+L_0^{-2+1/d+\eta}NB\}.\]
This will be of order $N^{1-\eta}$ if $\eta>0$ is sufficiently small,
and
\[\sup_{1\le t\le d/k}\rho(t)<1\;\;\;\mbox{and}\;\;\;\sup_{1\le t\le d/k} Q(t)<0\]
for
\[Q(t)=\left(-2+\frac{1}{d}\right)\frac{1+t-k(t)}{2}+t.\]
To handle the condition on $\rho(t)$ we note that the function 
attains its supremum at either $t=1$ ot $t=d/k$.  Moreover
if $v=k/d$ satisfies $5/9<v<1$ we find that $\rho(1)=9(1-v)/4<1$ and
\[\rho(d/k)=\frac{(1+v)(1-v^2)}{4v^2}.\]
This latter function is decreasing with respect to $v$, and takes the
value $196/225<1$ at $v=5/9$.  It follows that the supremum is
strictly less that 1 if $5/9<k/d<1$.

To verify the condition on $Q(t)$ we note that if $1\le t\le d/k$ then
\begin{eqnarray*}
Q'(t)&=&\left(-2+\frac{1}{d}\right)\frac{1}{2}\left\{1-
\frac{j}{2d}\left(1+\frac{d-kt}{j}+t\right)\left(-\frac{k}{j}+1\right)
\right\}+1\\
&=& \frac{1}{2d}-\left(1-\frac{1}{2d}\right)\frac{2k-d}{2d}
\left(1+\frac{d-kt}{j}+t\right)\\
&\le & \frac{1}{2d}-\left(1-\frac{1}{2d}\right)\frac{2k-d}{2d}
\left(1+\frac{d}{k}\right)\\
&<&0
\end{eqnarray*}
for $k>d/2$.  Thus 
\[Q(t)\le Q(1)=\frac{9j}{4d}\left(1-\frac{1}{2d}\right)
-\left(1-\frac{1}{d}\right),\]
which is strictly negative for 
\[j<\frac{4d}{9}\frac{2d-2}{2d-1}.\]
This condition is equivalent to
\[k>\frac{10d^2-d}{18d-9}=\frac{5d+2}{9}+\frac{2}{18d-9}.\]
Thus it is necessary and sufficient that
\[k\ge\frac{5d+3}{9}.\]
This completes the proof of Lemma \ref{mainlemma}, and hence also of
our two theorems.

Mathematical Institute,

24--29, St. Giles',

Oxford

OX1 3LB

UK
\bigskip

{\tt rhb@maths.ox.ac.uk}

\end{document}